\documentclass[10pt,draft]{amsart}
\usepackage{amsmath,amssymb,amsthm}
\begin{document}
\newtheorem{lem}{Lemma}[section]
\newtheorem{prop}{Proposition}[section]
\newtheorem{cor}{Corollary}[section]
\numberwithin{equation}{section}
\newtheorem{thm}{Theorem}[section]
\theoremstyle{remark}
\newtheorem{example}{Example}[section]
\newtheorem*{ack}{Acknowledgment}
\theoremstyle{definition}
\newtheorem{definition}{Definition}[section]
\theoremstyle{remark}
\newtheorem*{notation}{Notation}
\theoremstyle{remark}
\newtheorem{remark}{Remark}[section]
\newenvironment{Abstract}
{\begin{center}\textbf{\footnotesize{Abstract}}%
\end{center} \begin{quote}\begin{footnotesize}}
{\end{footnotesize}\end{quote}\bigskip}
\newenvironment{nome}
{\begin{center}\textbf{{}}%
\end{center} \begin{quote}\end{quote}\bigskip}

\newcommand{\triple}[1]{{|\!|\!|#1|\!|\!|}}
\newcommand{\xx}{\langle x\rangle}
\newcommand{\ep}{\varepsilon}
\newcommand{\al}{\alpha}
\newcommand{\be}{\beta}
\newcommand{\de}{\partial}
\newcommand{\la}{\lambda}
\newcommand{\La}{\Lambda}
\newcommand{\ga}{\gamma}
\newcommand{\del}{\delta}
\newcommand{\Del}{\Delta}
\newcommand{\sig}{\sigma}
\newcommand{\ome}{\omega}
\newcommand{\Ome}{\Omega}
\newcommand{\C}{{\mathbb C}}
\newcommand{\N}{{\mathbb N}}
\newcommand{\Z}{{\mathbb Z}}
\newcommand{\R}{{\mathbb R}}
\newcommand{\Rn}{{\mathbb R}^{n}}
\newcommand{\Rnu}{{\mathbb R}^{n+1}_{+}}
\newcommand{\Cn}{{\mathbb C}^{n}}
\newcommand{\spt}{\,\mathrm{supp}\,}
\newcommand{\Lin}{\mathcal{L}}
\newcommand{\SSS}{\mathcal{S}}
\newcommand{\F}{\mathcal{F}}
\newcommand{\xxi}{\langle\xi\rangle}
\newcommand{\eei}{\langle\eta\rangle}
\newcommand{\xei}{\langle\xi-\eta\rangle}
\newcommand{\yy}{\langle y\rangle}
\newcommand{\dint}{\int\!\!\int}
\newcommand{\hatp}{\widehat\psi}
\renewcommand{\Re}{\;\mathrm{Re}\;}
\renewcommand{\Im}{\;\mathrm{Im}\;}

\title[On the equipartition of energy]%
{On the equipartition of energy \\
for the critical $NLW$}

\author{}

\vspace{0.5cm}
\author[Luis Vega and Nicola Visciglia]{{Luis Vega\\
Universidad del Pais Vasco, Apdo. 64\\
48080 Bilbao, Spain\\
email: mtpvegol@lg.ehu.es\\
tel.:++34-946015475, fax: ++34-946012516\\
\vspace{0.3cm}
and\\
Nicola Visciglia\\
Dipartimento di Matematica Universit\`a di Pisa\\
Largo B. Pontecorvo 5, 56100 Pisa, Italy\\
email: viscigli@dm.unipi.it\\
tel.: ++39-0502212294, fax: ++39-0502213224}}

\maketitle
\date{}

Morawetz estimates,
dispersive estimates.

\vspace{0.1cm}

\begin{abstract}
We study some qualitative properties
of global solutions to the following focusing and defocusing critical $NLW$:
\begin{equation*}
\Box u+ \lambda u|u|^{2^*-2}=0, \hbox{ } \lambda\in {\mathbf R}
\end{equation*}
$$\hspace{2cm} u(0)=f\in \dot H^1( {\mathbf R}^n), 
\partial_t u(0)=g\in L^2( {\mathbf R}^n)$$
on ${\mathbf R}\times {\mathbf R}^n$ for $n\geq 3$,
where $2^*\equiv \frac{2n}{n-2}$.
We will consider the global solutions of the defocusing $NLW$ 
whose existence and scattering property is shown in \cite{shst}, \cite{sb}
and \cite{bg},
without any restriction on the initial data $(f,g)\in \dot H^1({\mathbf R}^n)
\times L^2({\mathbf R}^n)$.
As well as the solutions constructed in 
\cite{pecher} to the focusing $NLW$ for small initial data
and  to the ones obtained in
\cite{mk}, where a sharp condition on the smallness
of the initial data is given.
We prove that the solution $u(t, x)$ satisfies 
a family of  identities, that turn out to be a precised version
of the classical Morawetz estimates (see \cite{mor1}).
As a by--product we deduce that any global solution 
to critical $NLW$ belonging to a natural functional space satisfies:
$$\lim_{R\rightarrow \infty}\frac 1R \int_{\mathbf R} \int_{|x|<R} 
|\nabla_{x} u(t,x)|^2 \hbox{ }
dxdt
$$
$$=\lim_{R\rightarrow \infty}
\frac 1{2R} \int_{\mathbf R} \int_{|x|<R} (|\nabla_{t,x} u(t,x)|^2 
+ \frac{2 \lambda}{2^*} |u(t,x)|^{2^*}) \hbox{ } dxdt$$
$$=\int_{{\mathbf R}^n} (|\nabla_{t, x} u(0, x)|^2+
\frac{2 \lambda}{2^*} |u(0, x)|^{2^*}) \hbox{ } dx.$$
\end{abstract}

In this paper we study some qualitative properties of solutions to
the following family of Cauchy problems:
\begin{equation}\label{wave}
\Box u+ \lambda u|u|^{2^*-2}=0, \hbox{ } \lambda \in {\mathbf R}, 
\end{equation}
$$u(0)=f\in \dot H^1({\mathbf R}^n), \partial_t u(0)=g
\in L^2({\mathbf R}^n),$$
$$
(t, x) \in \mathbf R \times {\mathbf R}^n, n\geq 3,
$$
where $2^*\equiv \frac{2n}{n-2}$
and
$\Box\equiv \partial^2_{t}- \sum_{i=1}^n \partial_{x_i}^2.$

Notice that if $\lambda \equiv 0$, then \eqref{wave} 
reduces to the linear wave equation.
Since now on we shall refer to the Cauchy problem
\eqref{wave}
with $\lambda\geq 0$,  
as to the defocusing critical $NLW$
(similarly the focusing critical $NLW$ will be the Cauchy problem
\eqref{wave} with $\lambda<0$).

\vspace{0.1cm}

Along this paper we shall work with solutions $u(t, x)$ belonging
to the following space:
\begin{equation}\label{X}X\equiv {\mathcal C}({\mathbf R};\dot H^1({\mathbf R}^n))
\cap {\mathcal C}^1 ({\mathbf R}; L^2({\mathbf R}^n))
\end{equation}
\begin{equation*}
\cap L^\frac{2(n+1)}{n-2}({\mathbf R}\times {\mathbf R}^n)\cap L_{loc}^\frac{2(n+1)}{n-1}({\mathbf R};
\dot B^\frac 12_{\frac{2(n+1)}{n-1}}({\mathbf R}^n)).
\end{equation*}

We shall also assume that the conservation of the energy is satisfied
by the solutions $u(t, x)$, i.e.
\begin{equation}\label{fact}\int_{{\mathbf R}^n} (|\nabla_{t, x} u(T, x)|^2+
\frac{2 \lambda}{2^*} |u(T, x)|^{2^*}) \hbox{ } dx\equiv const
\hbox{ } \forall T\in \mathbf R.\end{equation}

\vspace{0.1cm}

Let us recall that the global well--posedness 
of the defocusing $NLW$ has been studied in \cite{mg}
provided that the initial data $(f, g)$ are regular.

Actually the global well--posedness of the defocusing Cauchy problem \eqref{wave}
has been studied in \cite{shst} for
initial data $(f, g)$ in the energy space 
$\dot H^1({\mathbf R}^n) \times L^2({\mathbf R}^n)$. 
In \cite{sb} and \cite{bg} 
the same problem has been analysed 
from the point of view of scattering theory (see also \cite{wo}). 
\noindent In particular by combining the results in \cite{sb} and
\cite{shst} 
it can be shown that for every $\lambda \geq 0$
and for every initial data $(f, g) \in \dot H^1({\mathbf R}^n) \times L^2({\mathbf R}^n)$,
there exists a unique solution
$u(t, x)$ to \eqref{wave} that belongs to the space
$X$ introduced in \eqref{X} and moreover \eqref{fact}
is satisfied.

\vspace{0.1cm}

Concerning the global well--posedness of the focusing $NLW$,
it is well--known that for every initial data
$(f, g)\in \dot H^1({\mathbf R}^n) \times L^2({\mathbf R}^n)$ small enough, i.e.
$$\int_{{\mathbf R}^n} 
(|\nabla_x f|^2 + |g|^2) \hbox{ } dx<\epsilon(|\lambda|)$$
for a suitable $\epsilon(|\lambda|)>0$, 
there exists a unique global solution to \eqref{wave} belonging to the space $X$ above
and moreover \eqref{fact} is satisfied.
For a proof of this fact see \cite{pecher}.

In the paper \cite{mk} a much more precised version
of the smallness assumption required on the initial data
is given in order to guarantee the global well--posedness 
of the focusing critical $NLW$.
In order to describe the result in \cite{mk}
let us introduce the function $W(x) \in \dot H^1({\mathbf R}^n)$
defined as follows:
$$W(x)\equiv \frac 1{( 1+\frac{|x|^2}{n(n-2)})^{\frac{n-2}2}}.
$$
Then in 
\cite{mk} it is shown that the  Cauchy problem  
\eqref{wave} with $\lambda=-1$, has a unique global solution
in the space $X$ introduced in \eqref{X},
provided that:
\begin{equation}\label{energymonot}
\int_{{\mathbf R}^n} (|\nabla_x f|^2 + |g|^2 - \frac{n-2}{n} |f|^{2^*}) \hbox{ }
dx < 
\int_{{\mathbf R}^n} (|\nabla_x W|^2 - \frac{n-2}{n} |W|^{2^*}) \hbox{ }
dx 
\end{equation}
and
\begin{equation}\label{energyprojecti}
\int |\nabla_x f|^2 \hbox{ } dx< \int |\nabla_x W|^2 \hbox{ } dx.
\end{equation}
Moreover in \cite{mk} it is proved that
blow--up occur provided that
$f$ and $g$ satisfy \eqref{energymonot}
and 
$\int |\nabla_x f|^2 \hbox{ } dx
> \int |\nabla_x W|^2 \hbox{ } dx$. 

It is also well--known that \eqref{fact} is satisfied by the solutions
constructed in \cite{mk}.

\vspace{0.1cm}

Our aim in this paper is to analyse some qualitative properties of global
solutions to \eqref{wave} in both focusing and defocusing case,
provided that such a global solutions exist and belong to the space $X$.
We are mainly interested
on the asymptotic behaviour for large $R>0$, of the following 
localized energies associated to the solutions $u(t, x)$ of \eqref{wave}:
\begin{equation}\label{moren}
\frac 1R \int_{\mathbf R} \int_{|x|<R} |\nabla_{x} u|^2 \hbox{ }
dxdt, \end{equation}
\begin{equation}\label{moren130}\frac 1R \int_{\mathbf R} \int_{|x|<R} |\nabla_{t,x} u|^2
\hbox{ }
dxdt
\hbox{ and }
\frac 1R \int_{\mathbf R} \int_{|x|<R} (|\nabla_{t,x} u|^2 + \frac{2\lambda}{2^*}|u|^{2^*})
\hbox{ }
dxdt.\end{equation}

Let us recall that the localized 
energies \eqref{moren} were first obtained  in 
\cite{kpv} following the ideas in \cite{AgmonHormmander}.
In this work we shall describe the asymptotic 
behaviour of the energies \eqref{moren} and \eqref{moren130}
as a consequence of a family of energy identities
satisfied by the global solutions to \eqref{wave}.

In our opinion those identities have its own 
interest since they represent a precised version
of the classical Morawetz inequalities, first proved
in \cite{mor1}.

\vspace{0.1cm}

Since now on we shall denote by $X$ the space 
defined in \eqref{X}.
Next we state the first result of this paper.

\begin{thm}\label{main}
Let $(f, g)\in \dot H^1({\mathbf R}^n)
\times L^2({\mathbf R}^n)$ and $\lambda \in \mathbf R$
be such that there exists a unique global solution 
$u(t, x)\in X$ to \eqref{wave}.
Assume moreover that $u(t, x)$ satisfies \eqref{fact}.
Let $\psi:{\mathbf R}^n\rightarrow \mathbf R$ be a 
radially symmetric function
such that:
$$(\sqrt{1+|x|^2}) \Delta_x \psi, \Delta^2_x \psi, 
\frac{\partial^2 \psi}{\partial x_i \partial x_j} 
 \in L^\infty({\mathbf R}^n) \hbox{ } \forall i, j=1,...,n$$ 
and
$$\lim_{|x|\rightarrow \infty}\partial_{|x|} \psi=\psi'(\infty).$$
Then we have
the following identity:
\begin{equation}\label{vai}\int_{\mathbf R}\int_{{\mathbf R}^n}
(\nabla_x u D^2_x \psi \nabla_x u - 
\frac 14 |u|^2 \Delta^2_x \psi  
+ \frac \lambda n|u|^{2^*}\Delta_x \psi  ) \hbox{ } dxdt
\end{equation}
$$=\psi'(\infty) \int_{{\mathbf R}^n} (|\nabla_x f|^2+
\frac{2 \lambda}{2^*} |f|^{2^*}+|g|^2) \hbox{ } dx.
$$
\end{thm}

\begin{remark} 
Let us point out that the hypothesis of theorem \ref{main}
are satisfied by the solutions constructed  in
\cite{shst} for defocusing $NLW$ and in \cite{pecher},\cite{mk}
for the focusing $NLW$.
\end{remark}

\begin{remark}
Let us underline that the identity \eqref{vai}
represents a precised version
of an inequality proved in \cite{mor1}, where 
\eqref{vai} is stated as an inequality 
and not as an identity in the special case $\psi\equiv |x|$.
\end{remark}

\begin{remark}
The same type of identities
as in theorem \ref{main}, have been proved in the context of the linear
Schr\"odinger equation in \cite{VV} and \cite{VVcmp}
respectively in the free and in the perturbative case.
The $L^2$--critical $NLS$ has been analysed from the same point of view
in \cite{VVindiana}. However the result stated for the critical $NLW$
in theorem \ref{main} is much more precise compared with the one
in \cite{VVindiana} for $NLS$. 

\noindent One of the main differences between $NLS$ and $NLW$
is that in the former case
an explicit representation of the asymptotic behavior of the solutions 
to the free Schr\"odinger equation is involved in the argument,
while in the case of $NLW$ it is not necessary. In this case one of the fundamental  
ingredients is the 
equipartition of energy, see  Propositon 2.1 below.
\end{remark}

\begin{remark}
Another difference between $NLW$ and $NLS$,
is that on the r.h.s.
in \eqref{vai} we get a quantity that is preserved along the evolution for $NLW$,
while in case of $NLS$ we get the 
$\dot H^\frac 12$--norm of the initial data, that
is not preserved along the evolution for $NLS$.
\end{remark}

Due to the freedom allowed
to the function $\psi$ in theorem \ref{main},
we can deduce the following result.

\begin{thm}\label{usu}
Let $(f, g)\in \dot H^1({\mathbf R}^n)
\times L^2({\mathbf R}^n)$, $\lambda \in \mathbf R$
and 
$u(t, x)\in X$ be as in theorem \ref{main}. Then
we have:
\begin{equation}\label{eqUCSL}\lim_{R\rightarrow \infty}
\frac 1R \int_{\mathbf R} \int_{|x|<R}
|\partial_{|x|} u|^2 \hbox{ } dxdt= \int_{{\mathbf R}^n} 
(|\nabla_x f|^2+
\frac{2 \lambda}{2^*} |f|^{2^*}+|g|^2) \hbox{ } dx.
\end{equation}
Moreover
\begin{equation}\label{eqUCSL123}\lim_{R\rightarrow \infty}
\frac 1R \int_{\mathbf R} \int_{|x|<R}
|\nabla_\tau u|^2 \hbox{ } dxdt= \lim_{R\rightarrow \infty}
\frac 1R \int_{\mathbf R} \int_{|x|<R}
|u|^{2^*} \hbox{ } dxdt=0,
\end{equation}
where $\partial_{|x|}$ and $\nabla_\tau$ represent
the radial derivative and the tangential part of the gradient respectively.
\end{thm}

Notice that theorem \ref{usu}
concerns mainly the concentration  of the spatial gradient of the solution.
Next we shall present another family of identities
that will allow us to study also the behaviour of the energies
connected with the time derivative of $u(t, x)$.

\begin{thm}\label{ut}
Let $(f, g)\in \dot H^1({\mathbf R}^n)
\times L^2({\mathbf R}^n)$, $\lambda \in \mathbf R$ and 
$u(t, x)\in X$ be as in theorem \ref{main}.
Let 
$\varphi: {\mathbf R}^n\rightarrow \mathbf R$ be a function
that satisfies the
following conditions:
$$\Delta_x \varphi, \langle x \rangle \varphi \in 
L^\infty({\mathbf R}^n).$$
Then the following identity holds
\begin{equation}\label{vaiut}\int_{\mathbf R} \int_{{\mathbf R}^n}
[( |\partial_t u|^2 - |\nabla_x u|^2
- \lambda |u|^{2^*})\varphi 
+ \frac 12 |u|^2 \Delta_x \varphi]  \hbox{ } dxdt=0.
\end{equation}
In particular we get:
\begin{equation}\label{equiptx}
\lim_{R\rightarrow \infty}
\frac 1R \int_{\mathbf R} \int_{|x|<R}
(|\nabla_x u|^2 -
|\partial_{t} u|^2) \hbox{ } dxdt=0,
\end{equation}
and 
\begin{equation}\label{vaiut2}\lim_{R\rightarrow \infty}
\frac 1R\int_{\mathbf R} \int_{|x|<R}
(|\nabla_{t,x} u|^2 
+ \frac{2\lambda}{2^*}|u|^{2^*}) \hbox{ } dxdt\end{equation}
$$=2 \int_{{\mathbf R}^n} (|\nabla_x f|^2+
\frac{2 \lambda}{2^*} |f|^{2^*}+|g|^2) \hbox{ } dx.
$$
\end{thm}

\begin{remark}
Notice that \eqref{equiptx} can be considered as a 
different version of the classical 
equipartition of the energy (see \cite{br}), whose classical version 
can be stated as follows: 
\begin{equation}\label{eqbr}\lim_{t\rightarrow \pm \infty} 
\int_{{\mathbf R}^n}
(|\partial_t u(t, x)|^2 - |\nabla_x u(t, x)|^2 )\hbox{ } dx =
0.\end{equation}
\end{remark}

\vspace{0.1cm}

\noindent Next we shall fix some notation.

\vspace{0.1cm}

{\bf Notation.}

\vspace{0.1cm}

\noindent For any $1\leq p, q\leq \infty$ 
$$L^p_x \hbox{ and } L^p_t L^q_x$$
denote the Banach spaces
$$L^p({\mathbf R}^n) \hbox{ and } L^p({\mathbf R}; L^q({\mathbf R}^n)),$$
and in the specific case $p=q$ we also use the notation
$$L^p_{t,x}\equiv  L^p({\mathbf R}; L^p({\mathbf R}^n)).$$
We shall denote by $L^{p,q}({\mathbf R}^n)$
the ususal Lorentz spaces  and by 
$\dot B^s_{p,2}({\mathbf R}^n)$
the Besov spaces.

\noindent For every $1\leq p\leq \infty$ we shall use the 
following mixed norm for functions defined on ${\mathbf R}^3$:
\begin{equation}\label{mixed}
\|f\|_{L^\infty_r L^p_\theta}^p 
\equiv \sup_{r>0} \int_{{\mathcal S}^2} |u(r\omega)|^p \hbox{ } d\omega
\end{equation}
where
$${\mathcal S}^2\equiv \{\omega\in {\mathbf R}^3| |\omega|=1\}$$
and $d\omega$ denotes the volume form on 
${\mathcal S}$.

\noindent We shall denote by
$\dot H^1_x$ the homogeneous Sobolev space $\dot H^1({\mathbf R}^n)$.

\noindent Given any couple of Banach spaces $Y$ and $Z$,
we shall denote by ${\mathcal L}(Y, Z)$
the space of linear and continuous
functionals between $Y$ and $Z$.
 
\noindent We shall denote by $${\mathcal C}_t(Y)
\hbox{ and }{\mathcal C}_t^1(Y)$$
respectively the spaces
$${\mathcal C}({\mathbf R};Y)
\hbox{ and } {\mathcal C}^1({\mathbf R};Y)$$
where $Y$ is a generic Banach space.

\noindent We shall denote by
$L^p_t(Y)$ the space
of $L^p$ functions defined on $\mathbf R$ and valued in $Y$.

\noindent We shall denote by $X$ the functional space
introduced in \eqref{X}.

\noindent Given a space--time dependent function $w(t, x)$ we shall denote by
$w(t_0)$ the trace of $w$ at fixed 
time $t\equiv t_0$, in case that it is 
well--defined.

\noindent We shall denote by $\int ...\hbox {  }dx ,
\int ... \hbox{ } dt$
and $\int \int ...\hbox{ } dx dt$ the 
integral of suitable functions
with respect to  space, time, 
and space--time variables respectively.

\noindent When it is not better specified we shall 
denote by $\nabla v$ the gradient
of any time--dependent function $v(t, x)$ 
with respect to the space variables.
Moreover $\nabla_\tau$ and $\partial_{|x|}$
shall denote respectively the angular 
gradient and the radial derivative.

\noindent If $\psi\in C^2({\mathbf R}^n)$, then 
$D^2 \psi$ will represent the hessian matrix of $\psi$.

\noindent Given a set $A\subset {\mathbf R}^n$ we denote by $\chi_A$
its characteristic function. 

\noindent The ball of radius $R>0$ in ${\mathbf R}^n$
shall be denoted as $B_R$.

\noindent We shall use the function
$$\langle x \rangle\equiv \sqrt{1+ |x|^2}.$$

\section{On the Strichartz estimates for critical $NLW$}

Recall that by combining the papers \cite{shst} and \cite{sb},
it follows that the defocusing $NLW$ is globally well--posed
in the Banach space introduced in \eqref{X} and moreover
the following properties hold: 
$$\lim_{t\rightarrow \pm \infty} \|u(t)\|_{L^{2^*}_x}=0$$
and
$$u(t, x)\in L^\frac{2(n+1)}{n-1}_t
\dot B^\frac 12_{\frac{2(n+1)}{n-1}}({\mathbf R}^n).$$
 
In the next proposition we gather some known facts that we shall use later on.
The main point is that
it applies to both focusing and defocusing $NLW$.

\begin{prop}\label{defocinterest}
Let $(f, g)\in \dot H^1({\mathbf R}^n)
\times L^2({\mathbf R}^n)$ and $\lambda \in \mathbf R$
be such that there exists a unique global solution 
$u(t, x)\in X$ to \eqref{wave}.
Then
we have:
\begin{equation}\label{conservenergy}
u(t,x)\in L^\infty_t \dot H^1_x;
\end{equation}
\begin{equation}\label{dispnlw}
\lim_{t\rightarrow \pm \infty} \|u(t)\|_{L^{2^*}_x}=0;
\end{equation}
\begin{equation}\label{strichgeneral}
u(t, x)\in L^\frac{2(n+1)}{n-1}_t
\dot B^\frac 12_{\frac{2(n+1)}{n-1}}({\mathbf R}^n).
\end{equation}
\end{prop}

\noindent{\bf Proof.}
For simplicity we shall prove \eqref{dispnlw} only in the case $t\rightarrow \infty$
and we shall also show boundedness of $\|u(t)\|_{\dot H^1_x}$
only for $t>0$.
The other cases can be treated in a similar way. 

\vspace{0.1cm}

{\em First step: $u(t, x)\in L^\infty_t \dot H^1_x$}

\vspace{0.1cm}

\noindent Since we are assuming $$u(t, x)\in X\subset L_{t,x}^{\frac{2(n+1)}{n-2}},$$
 we can deduce by standard techniques
in nonlinear scattering 
that $u(t,x)$ is asymptotically free. This means that
there exists $(f^+, g^+)\in \dot H^1_x \times L^2_x$
such that: 
\begin{equation}\label{train}
\lim_{t\rightarrow \infty} \|u(t)- u^+(t)\|_{\dot H^1_x}
+ \|\partial_t u(t) - \partial_t u^+(t)\|_{L^2_x}=0,
\end{equation}
where $$\Box u^+=0$$
$$u^+(0)=f^+, \partial_t u^+(0)=g^+.$$
The following computation is trivial:
\begin{equation}\label{train3}
\sup_{t\in \mathbf R}(\|\nabla_x u(t)\|_{L^2_x} 
+ \|\partial_t u(t)\|_{L^2_x})
\end{equation}
$$\leq
\sup_{t\in \mathbf R}\|\nabla_x u(t)- \nabla_x u^+(t)\|_{L^2_x} + 
\sup_{t\in \mathbf R}\|\nabla_x u^+(t)\|_{L^2_x}
$$
$$
+\sup_{t\in \mathbf R}\|\partial_t u(t)- \partial_t u^+(t)\|_{L^2_x} + 
\sup_{t\in \mathbf R}\|\partial_t u^+(t)\|_{L^2_x}
<\infty,$$
where at the last step we have used \eqref{train} and the conservation of the energy
for solutions to free wave equation.

\vspace{0.1cm}

{\em Second step: $u(t, x)\in L^\infty_t L^{2^*}_x$
and proof of \eqref{dispnlw}.}

\vspace{0.1cm}

\noindent By combining the previous step    
with the Sobolev embedding: 
\begin{equation}\label{sobolevcritico}
\dot H^1_x \subset L^{2^*}_x,\end{equation}
we deduce that 
$$u(t,x)\in L^\infty_t L^{2^*}_x.$$
On the other hand by combining 
again the Sobolev embedding with \eqref{train}
we get:
\begin{equation}\label{train2}
\lim_{t\rightarrow \infty}
\|u(t)\|_{L^{2^*}_x}
\leq \lim_{t\rightarrow \infty}
(\|u(t)-u^+(t)\|_{L^{2^*}_x} + \|u^+(t)\|_{L^{2^*}_x})=0, 
\end{equation}
where at the last step we have used proposition \ref{density}
in appendix I.

\vspace{0.1cm}

{\em Third step: proof of \eqref{strichgeneral}}

\vspace{0.1cm}

\noindent Once \eqref{dispnlw} has been shown, then 
the proof of \eqref{dispnlw} follows as in \cite{sb}.

\hfill$\Box$

\section{On the asymptotic behaviour of free waves}

First  we present a proposition whose content  
is well--known in the literature. However in  
Appendix II we shall present
a self--contained proof. 

\begin{prop}\label{dassios}
Let  $u(t,x)\in {\mathcal C}_t(\dot H^1_x)\cap 
{\mathcal C}^1_t(L^2_x)$ be the unique solution to:
$$\Box u=0 \hbox{ } (t, x) \in {\mathbf R} \times {\mathbf R}^n,
n\geq 3$$ 
$$u(0)=f\in \dot H^1_x, \partial_t u(0)=g\in L^2_x.$$
Then the following facts occur:
\begin{equation}\label{equip1gen}
\|\partial_t u(t) \pm \partial_{|x|} u(t)\|_{L^2_x}=o(1)
\hbox{ as } t\rightarrow \pm \infty,\end{equation}
\begin{equation}\label{equip2}
\lim_{t\rightarrow \pm \infty} \int |\nabla_x u(t)|^2 \hbox{ } dx
=\lim_{t\rightarrow \pm \infty} 
\int |\partial_t u(t)|^2 \hbox{ } dx
\end{equation}
\begin{equation*}
= \frac 12 \int (|\nabla_x f|^2 
+|g|^2) \hbox{ } dx;
\end{equation*}
\begin{equation}\label{equip3gen}
\int_{2|x|<|t|} (|\partial_t u(t)|^2 +|\nabla_x u(t)|^2) \hbox{ } dx =
o(1) \hbox{ as } t\rightarrow \pm \infty;
\end{equation}
\begin{equation}\label{equip37gen}
\int |\nabla_\tau u(t)|^2 \hbox{ } dx =
o(1) \hbox{ as } t\rightarrow \pm \infty;
\end{equation}
In particular, in the case $(f, g)\in C^\infty_0({\mathbf R}^n)
 \times C^\infty_0({\mathbf R}^n)$
we get the following stronger version of \eqref{equip1gen}
and \eqref{equip3gen}:
\begin{equation}\label{equip1}
\|\partial_t u(t) \pm  \partial_{|x|} u(t)\|_{L^2_x}=0\left( \frac 1
{|t|} \right)
\hbox{ as } t\rightarrow \pm \infty,\end{equation}
and
\begin{equation}\label{equip3}
\int_{2|x|<|t|} (|\partial_t u(t)|^2 +|\nabla_x u(t)|^2) \hbox{ } dx =
0\left ( \frac 1{t^2} \right) \hbox{ as } t\rightarrow \pm \infty.
\end{equation}
\end{prop}

\noindent{\bf Proof.} See Appendix II.

\hfill$\Box$

Next we shall study some asymptotic
expressions involving solutions to the free
wave equation with initial data in the energy
space $\dot H^1_x \times L^2_x$.
Those expressions will play a fundamental role in the sequel.
 
\begin{lem}\label{asym}
Assume that $u(t, x) \in {\mathcal C}_t (\dot H^1_x) \cap
{\mathcal C}_t^1(L^2_x)$ solves: 
\begin{equation}\label{linear}
\Box u=0
\end{equation}
$$u(0)=f\in \dot H^1_x, \partial_t u(0)=g\in L^2_x.$$

\noindent Let $\psi:{\mathbf R}^n \rightarrow \mathbf R$
be a radially symmetric function such that
the following limit exists:
$$\lim_{|x| \rightarrow \infty} \partial_{|x|}
\psi=\psi'(\infty)\in (0, \infty).$$
Then 
\begin{equation}
\label{gigl}\lim_{t\rightarrow \pm \infty}
\int \partial_t u(t) \nabla_x u(t) \cdot \nabla_x \psi \hbox{ } dx
= \mp \frac 12 \psi'(\infty) \int (|\nabla_x f|^2 + |g|^2) \hbox{ } dx.
\end{equation}
\end{lem}
\noindent{\bf Proof.}
We shall study only the case $t\rightarrow \infty$
(in fact the case $t\rightarrow -\infty$ reduces 
to the previous one since $v(t, x)\equiv 
u(-t, x)$ is still a solution to the free wave equation 
and its behaviour at infinity
is related to the behaviour of $u(t, x)$ as $t\rightarrow -\infty$).

\vspace{0.1cm}

\noindent Notice that we have:
\begin{equation}\label{fio}
\left |\int_{2|x|<t}  \partial_t u(t) \nabla_x u(t) 
\cdot \nabla_x \psi \hbox{ } dx\right |
\end{equation}
$$
\leq \frac 12 \|\nabla_x \psi\|_{L^\infty_x}
\int_{2|x|<t}
(|\partial_t u(t)|^2
+ |\nabla_x u(t)|^2) \hbox{ } dx$$
that due to 
\eqref{equip3gen} implies: 
$$\lim_{t\rightarrow \infty}
\int_{2|x|<t} \partial_t u(t)\nabla_x u(t) 
\cdot \nabla_x \psi \hbox{ } dx =0.
$$
Next notice that due to \eqref{equip1gen} we have
\begin{equation}\label{fio3}
\lim_{t\rightarrow \infty}
\int_{2|x|>t}  \partial_t u(t) 
\nabla_x u(t) \cdot \nabla_x \psi \hbox{ } dx
\end{equation}
$$= 
\lim_{t\rightarrow \infty}
\int_{2|x|>t}  \partial_t u(t) \partial_{|x|} u(t)  
\partial_{|x|}\psi \hbox{ } dx
= - \lim_{t\rightarrow \infty}
\int_{2|x|>t}  
|\partial_t u(t)|^2 \partial_{|x|}\psi \hbox{ } dx.
$$
On the other hand we have
$$\inf_{2|x|>t} 
(\partial_{|x|} \psi) \int_{2|x|>t} |\partial_t u(t)|^2 \hbox{ } dx
\leq 
\int_{2|x|>t}   
|\partial_t u(t)|^2 \partial_{|x|}\psi \hbox{ } dx
$$$$\leq 
\sup_{2|x|>t} 
(\partial_{|x|} \psi) 
\int_{2|x|>t} |\partial_t u(t)|^2 \hbox{ } dx.$$
Then due to the assumption done on $\partial_{|x|}\psi$ 
implies
\begin{equation}\label{fio4}\lim_{t\rightarrow \infty} 
\int_{2|x|>t}  
|\partial_t u(t)|^2 \partial_{|x|}\psi \hbox{ } dx
= \psi'(\infty) 
\lim_{t\rightarrow \infty} 
\int_{2|x|>t}  
|\partial_t u(t)|^2 \hbox{ } dx
\end{equation}
provided that the last limit exists.
By combining \eqref{fio}, \eqref{fio3}
and \eqref{fio4} we deduce that
the proof will be concluded provided that we can show 
\begin{equation}\label{fio5}\lim_{t\rightarrow \infty} 
\int_{2|x|>t}  
|\partial_t u(t)|^2 \hbox{ } dx= \frac 12 \int 
(|\nabla_x f |^2 + |g|^2) \hbox{ } dx.
\end{equation}
On the other hand we have
\begin{equation}\label{fio159}\lim_{t\rightarrow \infty} 
\int_{2|x|>t}  
|\partial_t u(t)|^2 \hbox{ } dx
\end{equation}
$$= \lim_{t\rightarrow \infty} 
\int  
|\partial_t u(t)|^2 \hbox{ } dx
- \lim_{t\rightarrow \infty} 
\int_{2|x|<t}  
|\partial_t u(t)|^2 \hbox{ } dx
$$
$$=\frac 12 \int 
(|\nabla_x f |^2 + |g|^2) \hbox{ } dx$$
where we have used \eqref{equip2} and \eqref{equip3gen}.

\hfill$\Box$

\begin{lem}\label{asym2}
Assume that $u(t, x) \in {\mathcal C}_t(\dot H^1_x)
\cap {\mathcal C}^1_t(L^2_x)$ solves \eqref{linear}.
Let $\phi:{\mathbf R}^n\rightarrow \mathbf R$
be a radially symmetric function such that:
$$\langle x \rangle \phi\in L^\infty_x.$$
Then we have
\begin{equation}\label{gigl2}
\lim_{t\rightarrow \pm \infty} \int \partial_t u(t) u(t) \phi \hbox{ }
dx=0.\end{equation}
\end{lem}

\noindent {\bf Proof.}
As in the proof of lemma \ref{asym}
it is sufficient to consider the limit for $t\rightarrow \infty$.

\vspace{0.1cm}

\noindent First notice that by combining
the decay assumption done on $\phi$
with the Hardy inequality we get
$$|\int \partial_t u(t) u(t) \phi \hbox{ } dx|
\leq \|\partial_t u(t)\|_{L^2_x} \|u(t)\phi\|_{L^2_x}
$$
$$\leq C (\|\partial_t u\|_{L^2_x}^2 + \|\nabla_x u\|_{L^2_x}^2) 
\equiv const \hbox{ } \forall t\in \mathbf R.$$
Due to this fact it is easy to show
that by a density argument 
it is sufficient to prove \eqref{gigl2}
under the assumptions
that $(f, g)\in C^\infty_0({\mathbf R}^n)
\times C^\infty_0({\mathbf R}^n)$.
Notice that if $u(t, x)$ is a regular solution to \eqref{linear},
then we have
$$\frac{d^2}{dt^2} \int |u(t)|^2 \hbox{ }
dx=2 \int (|\partial_t u(t)|^2 - |\nabla_x u(t)|^2) \hbox{ } dx,
$$
that due to \eqref{equip2} implies 
$$\frac{d^2}{dt^2} \int  |u(t)|^2 \hbox{ }
dx=o(1) \hbox{ as } t\rightarrow \infty.$$
After integration of this identity we get:
$$\int  |u(t)|^2 \hbox{ }
dx= \int  |f|^2 \hbox{ }
dx + 2 t (\int fg \hbox{ } dx) 
+ o(t^2),$$
and hence
\begin{equation}\label{ot2}\int  |u(t)|^2 \hbox{ }
dx=o(t^2) \hbox{ as } t\rightarrow \infty. 
\end{equation}
Next notice that we have: 
$$\int \partial_t u(t) u(t) \phi \hbox{ }
dx
=I(t) + II(t),
$$
where
$$
I(t)\equiv
\int_{2|x|>t} 
\partial_t u(t) u(t) \phi \hbox{ }
dx,$$
and
$$
II(t)\equiv
\int_{2|x|<t} 
\partial_t u(t) u(t) \phi \hbox{ }
dx.$$
Notice that due to the decay assumption done on $\phi$ 
we have
$$
I(t)\leq \frac C t \int_{2|x|>t} 
|u(t) \partial_t u(t)| \hbox{ } dx
$$
$$\leq \frac C t \|u(t)\|_{L^2_x} \|\partial_t u(t)\|_{L^2_x}
=o(1) \|\partial_t u(t)\|_{L^2_x}, 
$$
where we have used \eqref{ot2}. In particular we get: 
$$\lim_{t\rightarrow \infty} |I(t)|=0.$$
On the other hand
due to \eqref{ot2} and \eqref{equip3}
we have:
$$|II(t)|\leq \|\phi\|_{L^\infty_x}
\left(\int_{2|x|<t} |u(t)|^2 \hbox{ }dx\right)^\frac 12
\left(\int_{2|x|<t} |\partial_t u(t)|^2 \hbox{ }dx\right)^\frac 12
$$$$= Co(t)0(\frac 1t),$$
and hence
$$\lim_{t\rightarrow \infty} |II(t)|=0.$$
The proof is complete.

\hfill$\Box$

\section{Proof of theorem \ref{main}}

We shall need the following propositions.

\begin{prop}\label{nlw12}
Assume that $u(t, x)\in X$ is a global solution to \eqref{wave}
in dimension $n\geq 3$
for some $\lambda\in {\mathbf R}$ and $(f, g)\in \dot H^1_x\times L^2_x$.
Assume moreover that $u(t, x)$ satsifies \eqref{fact}.
Then we have:
$$\lim_{t\rightarrow \pm \infty}
\int \partial_t u(t) \nabla_x u(t)\cdot 
\nabla_x \psi(x)\hbox{ }
dx$$$$=\mp \frac 12 \psi'(\infty) \int_{{\mathbf R}^n}
(|\nabla_x f|^2 + \frac {2\lambda} {2^*} |f|^{2^*}+|g|^2 ) \hbox{ } dx,$$
where 
$\psi:{\mathbf R}^n \rightarrow \mathbf R$
is a radially symmetric function such that
the following limit exists
$$\lim_{|x| \rightarrow \infty} \partial_{|x|}
\psi=\psi'(\infty)\in (0, \infty),$$
\end{prop}

\noindent{\bf Proof.} 
As in the proof of proposition \ref{defocinterest} 
we only treat the case $t\rightarrow \infty$.
Let $u^+(t, x), f^+, g^+$ be as in the proof of proposition 
\ref{defocinterest}.

\noindent As a consequence of \eqref{train} we deduce that:
\begin{equation}\label{spagnolo}\lim_{t\rightarrow \infty}
\int \partial_t u(t) \nabla_x u(t)\cdot 
\nabla_x \psi \hbox{ }
dx\end{equation}
$$=\lim_{t\rightarrow \infty}
\int \partial_t u^+(t) \nabla_x u^+(t)\cdot 
\nabla_x \psi \hbox{ }
dx$$
$$=\lim_{t\rightarrow \infty}
-\frac 12 \psi'(\infty) \int (|\nabla_x f^+|^2 + |g^+|^2) \hbox{ } dx
$$
$$=\lim_{t\rightarrow \infty}
-\frac 12 \psi'(\infty) \int (|\nabla_x u^+(t)|^2 
+ |\partial_{t}u^+(t)|^2) \hbox{ } dx,
$$
where we have used
lemma \ref{asym} and the conservation of the energy
for the free wave equation.
Next notice that by combining \eqref{train}, \eqref{spagnolo},
the Sobolev embedding
$\dot H^1_x\subset L^{2^*}_x$ and the conservation of the energy for solutions to critical 
$NLW$ we get:
$$\lim_{t\rightarrow \infty}
\int \partial_t u(t) \nabla_x u(t)\cdot 
\nabla_x \psi \hbox{ }
dx$$
$$=\lim_{t\rightarrow \infty}
-\frac 12 \psi'(\infty) \int (|\nabla_x u^+(t)|^2 + |\partial_{t}u^+(t)|^2 
+ \frac{2\lambda}{2^*} |u^+(t)|^{2^*} - 
\frac{2\lambda}{2^*} |u^+(t)|^{2^*}) \hbox{ } dx$$
$$=\lim_{t\rightarrow \infty}
-\frac 12 \psi'(\infty) \int (|\nabla_x u(t)|^2 + |\partial_{t}u(t)|^2 
+ \frac {2\lambda}{2^*} |u(t)|^{2^*} - \frac{2\lambda}{2^*} |u^+(t)|^{2^*}) \hbox{ } dx$$
$$= -\frac 12 \psi'(\infty) \int (|\nabla_x f|^2 + |g|^2 
+ \frac {2\lambda}{2^*} |f|^{2^*}) \hbox{ } dx
+\frac \lambda {2^*} \psi'(\infty) \lim_{t\rightarrow \infty}
\int |u^+(t)|^{2^*}\hbox{ } dx$$
$$= -\frac 12 \psi'(\infty) \int (|\nabla_x f|^2 + |g|^2 
+ \frac{2\lambda}{2^*} |f|^{2^*}) \hbox{ } dx,$$
where we have used lemma \ref{gigl}
and the property
\begin{equation}\label{append2star}
\lim_{t\rightarrow \infty} \int |u^+(t)|^{2^*} \hbox{ } dx=0.
\end{equation}
The proof of \eqref{append2star} can be found in Appendix I.

\hfill$\Box$

\begin{prop}\label{nlw12345}
Assume that $u(t, x)\in X$ is a global solution to \eqref{wave}
in dimension $n\geq 3$
for some $\lambda\in {\mathbf R}$ and $(f, g)\in \dot H^1_x\times L^2_x$.
Let $\phi:{\mathbf R}^n\rightarrow \mathbf R$
be a radially symmetric function such that
$$\langle x \rangle \phi\in L^\infty_x,$$
then we have
\begin{equation}\label{gigl2345}
\lim_{t\rightarrow \pm \infty} \int \partial_t u(t) u(t) \phi \hbox{ }
dx=0.\end{equation}
\end{prop}
\noindent{\bf Proof.}
As in proposition \ref{defocinterest} 
it is sufficient to treat the case $t\rightarrow \infty$.

\vspace{0.1cm}

\noindent Let $u^+(t, x), f^+(x)$ and $g^+(x)$ be as in the proof of proposition
\ref{defocinterest}.
Notice that due to the decay assumption done on $\phi$ and 
due to the Hardy inequality we deduce that
\begin{equation}\label{giglietto}\lim_{t\rightarrow \infty} 
\|\phi(u(t)-u^+(t))\|_{L^2_x}
\end{equation}$$\leq 
C \lim_{t\rightarrow \infty} \|\nabla u(t) - \nabla u^+(t)\|_{L^2_x}=0$$
where at the last step we have used \eqref{train}.
On the other hand due again to \eqref{train} we have:
\begin{equation}\label{giglietto73}
\lim_{t\rightarrow \infty} \|\partial_t u(t) - \partial_t u^+(t)\|_{L^2_x}=0. 
\end{equation}
By combining \eqref{giglietto} and \eqref{giglietto73}
we deduce that
$$\lim_{t\rightarrow \infty} \int \partial_t u(t) u(t) \phi \hbox{ }
dx= \lim_{t\rightarrow \infty} \int \partial_t u^+(t) u^+(t) \phi \hbox{ }
dx=0$$
where at the last step we have used \eqref{gigl2}.

\hfill$\Box$

\noindent 
{\bf Proof of theorem \ref{main}.}
Following \cite{PV} we multiply the equation \eqref{wave}
by $\nabla_x \psi \cdot \nabla_x u + \frac 12 \Delta_x \psi u$
in order to get after integration by parts:
$$
\int_{-T}^T \int_{{\mathbf R}^n}
 \nabla_x u D^2_x \psi \nabla_x u 
- \frac 14 |u|^2 \Delta^2_x \psi  
+ \frac \lambda n |u|^{2^*}\Delta_x \psi \hbox{ } dxdt$$$$=
\sum_{\pm}
\left ( \mp \int_{{\mathbf R}^n} \partial_t u(\pm T) \nabla_x u(T) 
\cdot \nabla_x\psi 
+ \frac 12 \partial_t u(\pm T) u(\pm T)  \Delta_x \psi  \hbox{ } dx \right ).
$$
The proof can be completed by taking the limit as $T\rightarrow \infty$
and by using 
propositions \ref{nlw12} and \ref{nlw12345}.

\hfill$\Box$

\section{Proof of theorem \ref{usu}}

We start this section with some preliminary results
that will be useful
along the proof of theorem \ref{usu}.



\begin{prop}\label{hideo}
Let $(f, g)\in \dot H^1({\mathbf R}^n)
\times L^2({\mathbf R}^n)$ and $\lambda \in \mathbf R$
be such that there exists a unique global solution 
$u(t, x)\in X$ to \eqref{wave} for $n\geq 3$. 
Then:
\begin{equation}\label{radial}
\int\int \frac 1{\langle x \rangle} |u|^{2^*} \hbox{ }
dxdt<\infty,\end{equation}
and in particular:
\begin{equation}\label{utilsis}
\lim_{R\rightarrow \infty}
\frac 1R \int \int_{B_R} |u|^{2^*} \hbox{ } dxdt=0.
\end{equation}

\end{prop}

\noindent{\bf Proof.} 
Notice that due to the H\"older inequality in Lorentz spaces we get:
$$\int \frac 1{\langle x \rangle} |u|^{2^*} \hbox{ }
dx\leq \left \|\frac{1}{\langle x \rangle}\right \|_{L^{n, \infty}_x}
\|u\|_{L^{\frac{2^*n}{n-1},2^*}_x}^{2^*}\leq 
C \|u\|_{L^{\frac{2n(n+1)}{n^2-2n-1},2}({\mathbf R}^n)}^{2^*\theta}
\|u\|_{L^{2^*,2}({\mathbf R}^n)}^{2^*(1-\theta)},
$$
where:
$$\frac{\theta(n^2-2n-1)}{2n(n+1)}+\frac{1-\theta}{2^*}=\frac{n-1}{2^*n},$$
i.e. 
$$\theta=\frac{(n+1)(n-2)}{n(n-1)}.$$
By combining the previous inequality with the Sobolev embedding:
$$\dot H^1_x \subset L^{2^*,2}({\mathbf R}^n)$$
and
$$\dot B^\frac 12_{\frac{2(n+1)}{n-1},2}({\mathbf R}^n)
\subset L^{\frac{2n(n+1)}{n^2-2n-1},2}
({\mathbf R}^n),$$
we get:
$$\int \int \frac 1{\langle x \rangle} |u|^{2^*} \hbox{ }
dx\leq
C 
\|u\|_{L^\frac{2(n+1)}{n-1}_t\dot B^\frac 12_{\frac{2(n+1)}{n-1},2}({\mathbf R}^n)}^{\frac{2(n+1)}{n-1}}\|u\|_{L^\infty_t\dot H^1_x}^{\frac{4}{(n-1)(n-2)}}<\infty,$$
where at the last step we have used 
\eqref{conservenergy} and \eqref{strichgeneral}.
The proof of \eqref{radial} is complete.
Notice that \eqref{utilsis} follows by combining
\eqref{radial} with the dominated convergence theorem.

\vspace{0.1cm}





\hfill$\Box$




\begin{prop}\label{noma3}
Let $(f, g)\in \dot H^1_x\times L^2_x$ and $\lambda \in \mathbf R$
be such that there exists a unique global solution 
$u(t, x)\in X$ to \eqref{wave} with $n\geq 3$. 
Then we have:
\begin{equation}\label{due3}
\lim_{R\rightarrow \infty}  \int \int  
 |\Delta_x \phi_R| 
|u|^{2^*} \hbox{ } dxdt=0,
\end{equation}
where $\phi$ is a radially symmetric function
such that 
$$|\Delta_x \phi| \leq \frac C{\langle x \rangle} 
$$
and $\phi_R =R \phi\left (\frac{x}R\right )$.
\end{prop}

\noindent {\bf Proof.} 
By assumption we have
\begin{equation}\label{largevar} \int\int 
 |\Delta_x \phi_R ||u|^{2^*} \hbox{ } dxdt
\end{equation}
$$
\leq C \int
\int \frac C{R+ |x|}|u|^{2^*} \hbox{ } dxdt
\rightarrow 0 \hbox{ as }
R\rightarrow \infty,$$
where at the last step we have combined the dominated convergence theorem
with \eqref{radial}.

\hfill$\Box$

\begin{prop}\label{n=3}
Let $(f, g)\in \dot H^1_x\times L^2_x$ and $\lambda \in \mathbf R$
be such that there exists a unique global solution 
$u(t, x)\in X$ to \eqref{wave} for $n\geq 3$.
Then $u(t, x)$ satisfies:
\begin{equation}\label{sturom}\lim_{R\rightarrow \infty} \frac 1{R^3}
\int \int_{B_R} |u|^2 \hbox{ }
dxdt=0.\end{equation}
\end{prop}

In order to prove proposition \ref{n=3} 
we shall need some lemma.
In particular next result will be particularly useful
along the proof of proposition \ref{n=3} 
in the case $n=3$.

\begin{lem}\label{noc}
Let $u(t, x)\in {\mathcal C}_t(\dot H^1_x)\cap {\mathcal C}_t^1(L^2_x)$
be the unique solution to
$$\Box u=F\in L^1_tL^2_x, \hbox{ } (t, x)\in {\mathbf R}
\times {\mathbf R}^3$$
$$u(0)=f \in \dot H^1_x, \partial_t u(0)=g\in L^2_x.$$
For every $1\leq p <\infty$ there exists a constant $C\equiv C(p)>0$
such that
the following a--priori estimate holds:
$$\|u\|_{L^2_t L^\infty_r L^p_\theta}
\leq C (\|f\|_{\dot H^1_x} + \|g\|_{L^2_x}
+ \|F\|_{L^1_tL^2_x}).$$
\end{lem}

\noindent{\bf Proof.} 
In \cite{naoz} it is proved the following estimate for every $1\leq p<\infty$:
$$\|u\|_{L^2_t L^\infty_r L^p_\theta}\leq C 
(\|f\|_{\dot H^1_x} + \|g\|_{L^2_x})$$
where $u(t, x)\in {\mathcal C}_t (\dot H^1_x)\cap 
{\mathcal C}_t^1(L^2_x)$ is the unique solution to: 
$$\Box u=0$$
$$u(0)=f \in \dot H^1_x, \partial_t u(0)=g\in L^2_x.$$
The proof of lemma \ref{noc} in the case $F(t, x)\neq 0$ follows easily  
by combining the previous estimate with the Minkowski inequality 
and the Duhamel formula.

\hfill$\Box$

\begin{lem}\label{kenj}
Let $(f, g)\in \dot H^1_x\times L^2_x$ and $\lambda \in \mathbf R$
be such that there exists a unique global solution 
$u(t, x)\in X$ to \eqref{wave} for $n\geq 3$.
Then we have:
\begin{equation}\label{n>3}
\|u(t, x)\|_{L^2_t L_x^\frac{2n}{n-3}}<\infty
\hbox{ when } n\geq 4
\end{equation}
and 
\begin{equation}\label{n=3.}
\|u\|_{L^2_t L^\infty_r L^p_\theta}<\infty
\hbox{ } \forall 1\leq p<\infty \hbox{ when } n=3.
\end{equation}
\end{lem}

\noindent{\bf Proof.}
We split the proof in two parts.

\vspace{0.1cm}

{\em Proof of \eqref{n>3}}

\vspace{0.1cm}

\noindent Notice that by combining the Sobolev embedding:
$$\dot B^\frac 12_{\frac{2(n+1)}{n-1},2} 
\subset L^\frac{2n(n+1)}{n^2-2n-1}_x$$
with 
\eqref{strichgeneral}
we get
$$\|u(t, x)\|_{L^\frac{2(n+1)}{n-1}_t
L_x^{\frac{2n(n+1)}{n^2-2n-1}}}<\infty.$$
On the other hand
\eqref{conservenergy} implies
$$\|u(t, x)\|_{L^\infty_tL^{2^*}_x}<\infty$$
and hence 
by interpolation
\begin{equation}\label{chir}
\|u(t, x)\|_{L^\frac{2(n+2)}{n-2}_t L_x^\frac{2n(n+2)}{(n-2)(n+1)}}<\infty.
\end{equation}
Recall that $u(t, x)$ solves:
$$\Box u= - \lambda u|u|^{2^*-2}, \hbox{ } (t, x) 
\in {\mathbf R} \times {\mathbf R}^n, n\geq 4$$
$$u(0)=f\in \dot H^1({\mathbf R}^n), \partial_t u(0)=g
\in L^2({\mathbf R}^n),$$
and hence by Strichartz estimates (see \cite{gv})
we deduce:
\begin{equation}\label{STR}\|u(t, x)\|_{L^2_tL_x^{\frac{2n}{n-3}}} 
\leq C\left (\|f\|_{\dot H^1_x}
+\|g\|_{L^2_x}+ |\lambda | \|u^\frac{n+2}{n-2}\|_{L^2_tL_x^\frac{2n}{n+1}} \right).
\end{equation}
By combining this estimate with \eqref{chir} we get \eqref{n>3}.

\vspace{0.1cm}

{\em Proof of \eqref{n=3.}}

\vspace{0.1cm}

\noindent 
Notice that the proof of \eqref{n>3}
fails in dimension $n=3$ since in this case the end--point Strichartz estimate
(i.e. a version of \eqref{STR} for $n=3$) is false. 
Next we shall overcome this difficulty
by using lemma \ref{noc}.
By combining \eqref{strichgeneral} (where we choose $n=3$)
with the Sobolev embedding:
$$\dot B^\frac 12_{4,2}({\mathbf R}^3)
\subset L^{12}({\mathbf R}^3),$$
we deduce that 
\begin{equation}\label{atri}u(t, x)\in L^4_tL^{12}_x.\end{equation}
On the other hand due to \eqref{conservenergy} and due to the Sobolev embedding
$\dot H^1({\mathbf R}^3)\subset L^6({\mathbf R}^3))$,
we get $$\|u(t, x)\|_{L^\infty_t L^6_x}<\infty.$$
Hence by interpolation we get:
\begin{equation}\label{chir2}\|u\|_{L^5_tL^{10}_x}<\infty.\end{equation}
Next notice that $u(t, x)\in X$ solves the following Cauchy problem with forcing term:
$$\Box u= - \lambda u^5, \hbox{ } (t, x) 
\in {\mathbf R} \times {\mathbf R}^3$$
$$u(0)=f\in \dot H^1_x, \partial_t u(0)=g \in L^2_x.$$
By combining this fact with lemma \eqref{noc}
(where we choose $F=-\lambda u^5$) and \eqref{chir2}
we deduce:
$$\|u\|_{L^2_t L^\infty_r L^p_\theta}
\leq 
C(\|f\|_{\dot H^1_x} + \|g\|_{L^2_x}+
|\lambda| \|u\|_{L^5_tL^{10}_x}^5)<\infty.$$ 

\hfill$\Box$

\noindent{\bf Proof of proposition \ref{n=3}}
We shall prove proposition \ref{n=3} in dimension $n=3$
by using as a basic tool \eqref{n=3.}.
It will be clear that the same argument works in dimension $n>3$
provided that we use \eqref{n>3} instead of \eqref{n=3.}.

\vspace{0.1cm}

\noindent Since now on we shall assume $n=3$.
Notice that for every $T>0$ we have: 
$$\int_{T}^\infty \int_{B_R} |u|^2 dxdt=
\int_{T}^\infty  \int_{0}^R (\int_{{\mathcal S}^2}|u(t, r\omega)|^2  
d\omega )r^2 dr dt$$
$$\leq R^2 \int_{T}^\infty  
\int_{0}^R (\int_{{\mathcal S}^2}|u(t, r\omega)|^2 d\omega) drdt$$
$$\leq R^3 \int_{T}^\infty 
 \sup_{r\in (0,R)} \int_{{\mathcal S}^2}
|u(t, r\omega)|^2 d\omega dt $$
$$=R^3 \int (\sup_{r\in (0,\infty)} 
\|u(t, r\omega)\|_{L^2_\omega})^2
\hbox{ } dt=R^3 \|u\|_{L^2((T, \infty);L^\infty_r 
L^2_\omega)}^2.$$
By combining this fact with 
\eqref{n=3.} we get the following implication:
$$\forall \epsilon>0 \hbox{ there exists }
T_1(\epsilon)>0 \hbox{ s.t. }
\limsup_{R\rightarrow \infty}
\frac 1{R^3} \int_{T_1(\epsilon)}^\infty 
\int_{B_R} |u|^2 dxdt\leq \epsilon.
$$
Of course by a similar argument
we can prove that:
$$\forall \epsilon>0 \hbox{ there exists }
T_2(\epsilon)>0 \hbox{ s.t. }
\limsup_{R\rightarrow \infty}
\frac 1{R^3} \int^{-T_2(\epsilon)}_{-\infty} \int_{B_R} 
|u|^2 dxdt\leq \epsilon.
$$
In particular, if we choose $T(\epsilon)=\max\{T_1(\epsilon),
T_2(\epsilon)\}$, then we get:
\begin{equation}\label{maxmin}
\forall \epsilon>0 \hbox{ there exists }
T(\epsilon)>0 \hbox{ s.t. }
\end{equation}
\begin{equation*}
\limsup_{R\rightarrow \infty}
\frac 1{R^3} \int_{{\mathbf R}\setminus (-T(\epsilon); T(\epsilon)} 
\int_{B_R} |u|^2 dxdt\leq \epsilon.
\end{equation*}
Hence the proof of proposition \ref{n=3} (in the case $n=3$) will follow
from the following fact:
\begin{equation}\label{claimkenji}
\forall T>0 \hbox{ we have } \limsup_{R\rightarrow \infty} \frac 1{R^3}
\int_{-T}^T \int_{B_R}|u|^2 \hbox{ } dxdt=0.
\end{equation}
Notice that by using the H\"older inequality 
we get:
$$\int_{B_R}|u(t)|^2 \hbox{ } dx\leq R^2 \|u(t)\|_{L^6_x}^2,
$$
and this implies:
$$\frac 1{R^3}
\int_{-T}^T \int_{B_R}|u|^2 \hbox{ } dxdt
\leq \frac CR \int_{-T}^T \|u(t)\|_{L^6_x}^2 \hbox{ }
dt \leq \frac {2CT}R \|u\|_{L^\infty_t L^6_x}^2.$$
By combining this fact with \eqref{conservenergy} and with 
the Sobolev embedding
$\dot H^1_x\subset L^6_x$, we finally get  \eqref{claimkenji}.

\vspace{0.1cm}

\hfill$\Box$

\vspace{0.1cm}

\noindent {\bf Proof of theorem \ref{usu}}
First of all let us recall the following identity\begin{equation}\label{cube}\nabla_x \bar u 
D^2_x \psi \nabla_x u=
\partial_{|x|}^2\psi |\partial_{|x|} u|^2 + \frac{\partial_{|x|} 
\psi}{|x|}|\nabla_\tau u|^2,\end{equation}
where $\psi$ is a radially symmetric function.
By using this identity and by choosing 
in the identity \eqref{vai}
the function $\psi\equiv \langle x \rangle$, then it is easy to deduce that
\begin{equation*}
\int \int_{|x|>1} \frac{|\nabla_\tau u|^2}{|x|} \hbox{ } dxdt <\infty.
\end{equation*}
In particular we deduce:
\begin{equation}\label{revised}
\lim_{R\rightarrow \infty}\int \int_{|x|>R} \frac{|\nabla_\tau u|^2}{|x|} \hbox{ } dx=0.
\end{equation}
and
\begin{equation}\label{ultimissim1792}
\lim_{R\rightarrow \infty}
\frac 1R \int \int_{B_R} |\nabla_\tau u|^2 \hbox{ } dxdt=0.
\end{equation}
By combining \eqref{utilsis} with \eqref{ultimissim1792}
we get \eqref{eqUCSL123}.

\vspace{0.1cm}

\noindent Next we shall prove
\eqref{eqUCSL}.
For any $k\in \mathbf N$
we fix a function
$h_k(r)\in C^\infty_0(\mathbf R; [0, 1])$ such that:
\begin{equation}\label{hk}h_k(r)=1 \hbox{  } \forall r\in {\mathbf R}
\hbox{ s.t. } |r|<1, h_k(r)=0 \hbox{  } \forall r \in {\mathbf R}
\hbox{ s.t. } |r|>\frac{k+1}{k}, 
\end{equation}
$$h_k(r)= h_k(-r) \hbox{  } \forall r \in \mathbf R.$$
Let us introduce the functions
$\psi_k(r), H_k(r) \in C^\infty( {\mathbf R})$:
\begin{equation}\label{kr}\psi_k(r) = \int_0^r (r-s) h_k(s) ds \hbox{  } \hbox{ and } \hbox{  }
H_k(r)=\int_0^r h_k(s) ds.\end{equation}
Notice that 
\begin{equation}\label{differ}
\psi_k''(r)=h_k(r), \psi_k'(r)= H_k(r) \forall r\in {\mathbf R} \hbox{ and }
\lim_{r\rightarrow \infty} \partial_r \psi_k(r)=\int_0^\infty h_k(s)ds.
\end{equation}
Moreover an elementary computation 
shows that: 
\begin{equation*}
\Delta_x \psi_k\leq \frac C{\langle x \rangle} 
\hbox{ } \forall x\in 
{\mathbf R}^n , n\geq 3 
\end{equation*}
and
\begin{equation}\label{4dref}
\Delta^2_x \psi_k(x)= 
\frac C{|x|^3} \hbox{  } \forall x\in 
{\mathbf R}^n \hbox{ s.t. } |x| \geq 2 \hbox{ and } n\geq 4,
\end{equation}
\begin{equation}\label{3dref}
\Delta^2_x \psi_k(x)= 0 \hbox{  } \forall x\in 
{\mathbf R}^3 \hbox{ s.t. } |x| \geq 2,
\end{equation}
where $\Delta^2_x$ is the bilaplacian operator.
Thus the functions
$\phi\equiv \psi_k$ satisfy the assumptions of proposition
\ref{noma3}.
In the sequel we shall need the rescaled functions
\begin{equation}\label{kR}\psi_{k, R}(x)
\equiv R \psi_k\left (\frac{x}{R} \right ) \forall x\in {\mathbf R}^n, k\in 
{\mathbf N} \hbox{ and }R>0,
\end{equation}
where $\psi_k$ is defined in \eqref{kr}.
Notice that by combining the general identity 
\eqref{cube} with \eqref{vai}, where we choose $\psi= 
\psi_{k,R}$ defined in \eqref{kR}, and recalling \eqref{differ} we get:
\begin{equation}\label{limit}\int \int
\left ( \partial_{|x|}^2 \psi_{k,R} |\partial_{|x|} u|^2 + 
\frac{\partial_{|x|} \psi_{k,R}}
{|x|} | \nabla_\tau u|^2
\right .\end{equation}
\begin{equation*}
\left . - \frac 14
|u|^2 \Delta^2_x \psi_{k,R} + 
\frac \lambda n |u|^{2^*}\Delta_x \psi_{k,R} \right )dx dt 
\end{equation*}
\begin{equation*}
= \left (\int_0^\infty h_{k}(s) \hbox{ } ds
\right )\int (|\nabla_x f|^2+
\frac{2 \lambda}{2^*} |f|^{2^*}+|g|^2) \hbox{ } dx 
\hbox{ } \forall k\in {\mathbf N}, R>0.
\end{equation*}
Notice also that due to \eqref{3dref}
we get:
$$ \int\int_{{\mathbf R}^3}
|\Delta^2_x \psi_{k,R}|  |u|^2 \hbox{ } dxdt
\leq \frac C{R^3}\int\int_{B_R} |u|^2 \hbox{ } dxdt
$$
provided that $n=3$, and in particular by using
\eqref{sturom} we get
\begin{equation}\label{coltar3}\lim_{R\rightarrow \infty}
\int\int_{{\mathbf R}^3}
|\Delta^2_x \psi_{k,R}|  |u|^2 \hbox{ } dxdt=0.\end{equation}
In the case $n\geq 4$ we use \eqref{4dref}
in order to deduce:
\begin{equation}\label{matmai} \int\int_{{\mathbf R}^n}
|\Delta^2_x \psi_{k,R}|  |u|^2 \hbox{ } dxdt
\end{equation}$$\leq C \left (\frac 1{R^3}\int\int_{B_R} |u|^2 \hbox{ } dxdt
+ \int\int_{{\mathbf R}^n\setminus B_R}
 \frac{|u|^2}{|x|^3} \hbox{ } dxdt \right ).$$
On the other hand an explicit computation
shows that if we choose in \eqref{vai}
$\psi\equiv \langle x\rangle$, when $n\geq 4$,
then we get:
\begin{equation}\label{origmor}
\int \int_{{\mathbf R}^n} \frac{|u|^2}{|x|^3} \hbox{ }
dxdt<\infty \hbox{ for } n\geq 4,\end{equation}
that in conjunction with the Lebesgue dominated convergence theorem,
\eqref{sturom} and \eqref{matmai} implies:
\begin{equation}\label{coltar4}\lim_{R\rightarrow \infty} \int\int_{{\mathbf R}^n}
|\Delta^2_x \psi_{k,R}|  |u|^2 \hbox{ } dxdt=0
\hbox{ for } n\geq 4.\end{equation}
By using \eqref{coltar3}, \eqref{coltar4}, \eqref{due3}
and \eqref{revised} we get:
\begin{equation}\label{modifn=3}
\lim_{R\rightarrow \infty} \int\int
\left (\partial_{|x|} \psi_{k,R}\frac{|\nabla_\tau u|^2}{|x|} 
- \frac 14 \Delta^2_x \psi_{k,R}  |u|^{2} 
+ \frac \lambda n \Delta_x \psi_{k,R}
|u|^{2^*} \right ) \hbox{ } dx dt=0
\end{equation}
for every $k\in {\mathbf N}$ and for every dimension $n\geq 3$.
We can combine this fact with \eqref{limit} in order to deduce:
\begin{equation}\label{crusia}\lim_{R\rightarrow \infty} \int\int
\partial_{|x|}^2\psi_{k,R} |\partial_{|x|} u|^2
dx dt \end{equation}
\begin{equation*}=\left (\int_0^\infty h_{k}(s) ds
\right )\int (|\nabla_x f|^2+
\frac{2 \lambda}{2^*} |f|^{2^*}+|g|^2) 
\hbox{ } dx \hbox{ } \forall k\in {\mathbf N}.
\end{equation*}
On the other hand,
due to the properties of $h_{k}$ (see \eqref{hk}), we get
$$\frac 1R \int \int_{B_R}|\partial_{|x|} u|^2
dx dt \leq \int \int
\partial_{|x|}^2 \psi_{k,R}  |\partial_{|x|} u|^2  dt dx $$$$
= \frac 1R \int \int
h_k\left( \frac {x}R\right)|\partial_{|x|} u|^2  dt dx 
\leq \frac 1R \int\int_{|x|<\frac {k+1}{k}R} |\partial_{|x|} u|^2 dxdt$$
that due to \eqref{crusia} implies:
\begin{equation}\label{limsup}\limsup_{R\rightarrow \infty} \frac 1R \int 
\int_{B_R}|\partial_{|x|} u|^2 dxdt\end{equation}
$$ 
\leq  \left (\int_0^\infty h_k(s) ds \right ) \int 
(|\nabla_x f|^2+
\frac{2 \lambda}{2^*} |f|^{2^*}+|g|^2) \hbox{ } dx
$$$$\leq\frac {k+1}{k}
\liminf_{R\rightarrow \infty} \frac 1R \int\int_{B_R}|\partial_{|x|} u|^2 dxdt
\hbox{  } \forall k \in {\mathbf N}.$$
Since $k\in \mathbf N$ is arbitrary and since 
the following identity is trivially satisfied:
$$\lim_{k\rightarrow \infty}\int_0^\infty h_k(s) ds=1,$$
we can deduce 
\eqref{eqUCSL} by using \eqref{limsup}.

\vspace{0.1cm}

\noindent The proof is complete.

\hfill$\Box$

\section{Proof of theorem \ref{ut}}

\vspace{0.1cm}

{\em First step: proof of \eqref{vaiut}}

\vspace{0.1cm}

\noindent Following \cite{PV} we multiply the equation \eqref{wave}
by $\varphi u$ and integrating the corresponding identity
on the strip $(-T, T)$ we get:
\begin{equation}\label{vaiut123456}\int_{-T}^T 
\int_{{\mathbf R}^n}
( |\partial_t u|^2 - |\nabla_x u|^2
- \lambda |u|^{2^*})\varphi 
+ \frac 12 |u|^2 \Delta_x \varphi  \hbox{ } dxdt
\end{equation}\begin{equation*}
=\sum_{\pm} \pm \int \partial_t u(\pm T) u(\pm T) \varphi \hbox{ } dx 
\end{equation*}
By taking the limit as $T\rightarrow \infty$
and by using proposition \ref{nlw12345} we get 
\eqref{vaiut}.

\vspace{0.1cm}

{\em Second step: proof of \eqref{equiptx}}

\vspace{0.1cm}

\noindent For any $k\in \mathbf N$
we fix a function
$\varphi_k(r)\in C^\infty_0(\mathbf R; [0, 1])$ such that:
\begin{equation}\label{varphik}\varphi_k(r)=1 \hbox{  } \forall r\in {\mathbf R}
\hbox{ s.t. } |r|<1, \varphi_k(r)=0 \hbox{  } \forall r \in {\mathbf R}
\hbox{ s.t. } |r|>\frac{k+1}{k}, 
\end{equation}
$$\varphi_k(r)= \varphi_k(-r) \hbox{  } \forall r \in \mathbf R.$$
We also introduce the rescaled functions
$$\varphi_{k,R}\equiv \frac 1R \varphi_{k}\left(\frac xR \right).$$
Notice by combining the cut--off property of the functions $\varphi_k$
with \eqref{utilsis}
and \eqref{sturom}, we get:
$$\lim_{R\rightarrow \infty}
\int \int |u|^{2^*}\varphi_{k,R} \hbox{ } dxdt=\lim_{R\rightarrow \infty}
\int \int |u|^2 \Delta_x \varphi_{k,R} \hbox{ } dxdt=0 \hbox{ } \forall k\in \mathbf N$$
in any dimension $n\geq 3$.
By using this fact in conjunction with \eqref{vaiut}, 
where we choose $\varphi\equiv \varphi_{k,R}$,
we get:
\begin{equation}\label{confed}
\lim_{R\rightarrow \infty} \int \int
(|\partial_t u|^2 - |\nabla_x u|^2)
\varphi_{k,R} \hbox{ } dxdt=0 \hbox{ } \forall k \in \mathbf N. 
\end{equation}
Notice that by combining \eqref{confed}
with the cut--off properties of $\varphi_k$ we get:
\begin{equation}\label{stama}
\forall k\in {\mathbf N} \hbox{ there exists } R(k)>0 \hbox{ s.t. }
\end{equation}
$$\frac 1R \int\int_{B_R}|\partial_t u|^2 \hbox{ } dxdt
\leq \frac {k+1}k \frac{1}{R(\frac{k+1}{k})}
\int \int_{B_{R(\frac{k+1}{k})}} |\nabla_x u|^2 \hbox{ } dxdt
+ \frac 1k \hbox{ } \forall R>R(k).$$
By combining \eqref{stama} with \eqref{eqUCSL} and \eqref{eqUCSL123} ,
we get:
$$\limsup_{R\rightarrow \infty}
\frac 1R \int\int_{B_R}|\partial_t u|^2 \hbox{ } dxdt$$$$\leq 
\frac{k+1}{k} \int 
(|\nabla_x f|^2+
\frac{2 \lambda}{2^*} |f|^{2^*}+|g|^2) \hbox{ } dx
+ \frac 1k \hbox{ } \forall k\in \mathbf N$$
and in particular
$$\limsup_{R\rightarrow \infty}
\frac 1R \int\int_{B_R}|\partial_t u|^2 \hbox{ } dxdt$$$$\leq 
\int 
(|\nabla_x f|^2+
\frac{2 \lambda}{2^*} |f|^{2^*}+|g|^2) \hbox{ } dx.$$
Similarly one can show that:
$$\liminf_{R\rightarrow \infty}
\frac 1R \int\int_{B_R}|\partial_t u|^2 \hbox{ } dxdt$$$$\geq 
\int_{{\mathbf R}^n} 
(|\nabla_x f|^2+
\frac{2 \lambda}{2^*} |f|^{2^*}+|g|^2) \hbox{ } dx,$$
and finally we get:
$$\lim_{R\rightarrow \infty}
\frac 1R \int\int_{B_R}|\partial_t u|^2 \hbox{ } dxdt= 
\int
(|\nabla_x f|^2+
\frac{2 \lambda}{2^*} |f|^{2^*}+|g|^2) \hbox{ } dx
$$
$$= \lim_{R\rightarrow \infty}
\frac 1R \int\int_{B_R}|\nabla_x u|^2 \hbox{ } dxdt$$
where at the last step we have combined \eqref{eqUCSL}
with \eqref{eqUCSL123}.
The proof of 
\eqref{equiptx} is complete.

\vspace{0.1cm}

\noindent Finally notice that by combining \eqref{eqUCSL}, \eqref{eqUCSL123}
and \eqref{equiptx}, we get \eqref{vaiut2}.

\hfill$\Box$

\section{Appendix I}

The aim of this appendix is to show that
the $L^{2^*}$-norm of the solution to the following Cauchy problem: 
\begin{equation}\label{lwfg}\Box u=0
\end{equation}
$$u(0)=f \in \dot H^1_x, \partial_t u(0)=g\in L^2_x,$$
goes to zero as $t\rightarrow \pm \infty$.
Notice that this fact represents
a slight improvement compared
with the usual Strichartz estimate
$$\|u(t, x)\|_{L^\infty_t L^{2^*}_x} \leq C 
(\|f\|_{\dot H^1_x} + \|g\|_{L^2_x}).$$
On the other hand in proposition \ref{noratedecay} 
we shall show that in general no better
result can be expected. In fact we shall show that there cannot exist 
a--priori any rate on the decay  of the
$L^{2^*}$-norm of the solution to \eqref{lwfg}.

\vspace{0.1cm}

\noindent Along this section, when it is not better specified,
we shall denote by 
$T(t)(f, g)$ the solution to the Cauchy problem \eqref{lwfg} 
with initial data $(f, g)$ computed at time $t$, i.e.:
$$T(t): \dot H^1_x\times L^2_x \ni (f, g) \rightarrow u(t)\in \dot H^1_x,$$
where $u(t, x)$ solves \eqref{lwfg}.

\begin{prop}\label{density}
Let $u(t, x)\in {\mathcal C}_t (\dot H^1_x)\cap 
{\mathcal C}^1_t (L^2_x)$ be the unique solution to \eqref{lwfg},
then
$$\lim_{t\rightarrow \pm \infty} \|u(t)\|_{L^{2^*}_x}=0.$$
\end{prop}

\noindent{\bf Proof.}
We treat for simplicity the case $t\rightarrow \infty$
(the case $t\rightarrow -\infty$ can be treated in a similar way).
Notice that due to the Sobolev embedding and 
the conservation of the energy we have:
\begin{equation}\label{uniform}
\|u(t)\|_{L^{2^*}_x}^2
\leq S (\|\nabla_x u(t)\|_{L^2_x}^2
+ \|\partial_t u(t)\|_{L^2_x}^2)
\end{equation}\begin{equation*}
= S (\|\nabla_x f\|_{L^2_x}^2 + \|g\|_{L^2_x}^2)
\hbox{ } \forall (f, g)\in \dot H^1_x \times L^2_x
\end{equation*}
In particular the operators
$T(t)$ introduced above, are 
uniformly bounded for every $t>0$
in the space ${\mathcal L}(\dot H^1_x \times L^2_x, L^{2^*}_x)$.

\noindent On the other hand we have the following dispersive estimate (see \cite{SS}):
\begin{equation}\label{shatastruwe}\|u(t, x)\|_{L^\infty_x}
\leq \frac C{t^{\frac{n-1}2}}
(\|f\|_{\dot B_{1,1}^{\frac{m-1}2}}
+ \|g\|_{\dot B_{1,1}^{\frac{m+1}2}})\end{equation}
(here $\dot B^{s}_{p,q}$ denotes the
standard Besov spaces). Notice also that
the Fourier representation of the solution to \eqref{lwfg} 
implies:
\begin{equation}\label{fourier}
\|u(t)\|_{L^2_x}\leq C (\|f\|_{L^2_x} + \|g\|_{\dot H^{-1}_x}).
\end{equation}
In particular by combining
\eqref{shatastruwe} with \eqref{fourier}
we deduce:
$$\|u(t)\|_{L^{2^*}_x} \leq \|u(t)\|_{L^\infty_x}^{\frac 2n}
\|u\|_{L^2_x}^{\frac {n-2}n}
\leq \frac{C}{t^\frac{n-1}n} \hbox{ }
\forall (f, g)\in C^\infty_0({\mathbf R}^n) \times C^\infty_0({\mathbf R}^n),$$
where $C\equiv C(f, g)>0$.
As  a consequence we get 
\begin{equation}\label{dense}\lim_{t\rightarrow \infty}
\|u(t)\|_{L^{2^*}}=0  \hbox{ }
\forall (f, g)\in C^\infty_0({\mathbf R}^n) \times C^\infty_0({\mathbf R}^n).
\end{equation}
It is now easy to remove in \eqref{dense}
the regularity assumption
$(f, g)\in  C^\infty_0({\mathbf R}^n) \times C^\infty_0({\mathbf R}^n)$ by a 
classical density argument.

\hfill$\Box$

\vspace{0.1cm}

Notice that the previous result represents
a slight improvement compared
with the usual Strichartz estimate:
$$\|u(t, x)\|_{L^\infty_t L^{2^*}_x} \leq C 
(\|f\|_{\dot H^1_x} + \|g\|_{L^2_x}).$$

On the other hand next proposition shows that in general no better
result can be expected, since there cannot exist 
a--priori any rate on the decay  of the
$L^{2^*}$-norm of the solution to \eqref{lwfg}.

\begin{prop}\label{noratedecay}
Let $\gamma\in {\mathcal C} ([0, \infty); {\mathbf R})$
be any function such that
$$\lim_{t\rightarrow \infty}\gamma (t)=\infty.$$
Then there exists $g\in L^2_x$ such that:
$$\|u(t_n)\|_{L^{2^*}_x}>\frac  1{\gamma(t_n)},$$
where $\{t_n\}$ is a suitable sequence that goes to $+\infty$
and
\begin{equation}\label{CP}
\Box u=0\end{equation} $$u(0)=0, \partial_t u (0)=g\in L^2_x.$$
\end{prop}

\noindent{\bf Proof.}
We claim the following fact:

\begin{equation}\label{claim} 
\|S(t)\|_{{\mathcal L}(L^2_x, L^{2^*}_x)}\geq \epsilon_0>0 \
\hbox{ where }
\end{equation}$$S(t): L^2_x\ni g \rightarrow u(t)\in L^{2^*}_x
$$$$\hbox{ is the solution operator associated to the Cauchy problem
\eqref{CP}}.
$$
Notice that due to \eqref{claim} we get:
$$\lim_{t\rightarrow \infty}
\gamma(t)\|S(t)\|_{{\mathcal L}(L^2_x, L^{2^*}_x)}=\infty,$$
and in particular
due to the Banach--Steinhaus theorem the operators
$\gamma(t) S(t)$ cannot be pointwisely bounded,
or in an equivalent way there exists at least one $g\in L^2_x$ such that:
\begin{equation}\label{BS}\sup_{[0, \infty)} \gamma(t)\|S(t)g\|_{L^{2^*}_x}=\infty.
\end{equation}
On the other hand the function $\gamma(t)\|S(t)g\|_{L^{2^*}_x}$
is bounded on bounded sets of $[0, \infty)$,
and hence \eqref{BS} implies that
$$\limsup_{t\rightarrow \infty} \gamma(t)\|S(t)g\|_{L^{2^*}_x}=\infty$$
and it completes the proof.

\vspace{0.1cm}

\noindent Next we shall prove \eqref{claim}.
Let us fix $h\in L^2_x$ such that:
$$\|h\|_{L^2_x}=1 \hbox{ and } \|S(1)h\|_{L^{2^*}_x}
=\|u(1, x)\|_{L^{2^*}_x}
=\eta_0>0$$
where $u(t, x)$ denotes the unique solution
to \eqref{CP} with $g=h$. 

\vspace{0.1cm}

\noindent A rescaling argument implies
that 
$u_{\epsilon}(t, x)\equiv \epsilon^{\frac n2-1} 
u(\epsilon t, \epsilon x)$ solves \eqref{CP}
with initial data $g\equiv h_\epsilon\equiv \epsilon^{\frac n2} h(\epsilon x)$.
In particular this implies that:
$$S\left (\frac 1 \epsilon\right ) h_\epsilon=
\epsilon^{\frac n2-1}u(1, \epsilon x) \hbox{ and }
\|h_\epsilon\|_{L^2_x}=1,
$$
and hence:
$$\|S\left (\frac 1 \epsilon\right )\|_{{\mathcal L}(L^2_x, L^{2^*}_x)}
\geq \|S\left (\frac 1 \epsilon\right ) h_\epsilon\|_{L^{2^*}_x}
$$$$=\epsilon^{\frac n2-1}\| u(1, \epsilon x) 
\|_{L^{2^*}_x}=
\|u(1)\|_{L^{2^*}_x}=\|S(1)h\|_{L^{2^*}_x}=\eta_0>0 \hbox{ } \forall \epsilon>0.$$
The proof of \eqref{claim} is complete.

\hfill$\Box$

\section{Appendix II}

This section is devoted to the proof of proposition
\ref{dassios}. Let us underline that its content 
is well--known in the literature,
in particular it
contains the equipartition of the energy principle
first proved in \cite{br} by using Fourier
analysis.

\vspace{0.1cm}

The aim of this section is to present a proof
that involves the conformal energy. 

\vspace{0.1cm}

\noindent{\bf Proof of prop. \ref{dassios}}
First of all notice that \eqref{equip1gen}
implies:
$$\lim_{t\rightarrow \infty} \int |\partial_t u(t)|^2
= \lim_{t\rightarrow \infty} \int |\partial_{|x|} u(t)|^2 \hbox{ } dx.$$
By combining this fact with \eqref{equip2}
and with the following trivial identity: 
$$|\nabla_x u|^2= |\partial_{|x|}u|^2 + |\nabla_\tau u|^2,$$
we can deduce \eqref{equip37gen}.
Hence it is enough to prove \eqref{equip1gen} and \eqref{equip2}
in order to deduce \eqref{equip37gen}.

\vspace{0.1cm}

\noindent Next notice that by a density argument it is sufficient to prove
\eqref{equip2}, \eqref{equip1} and \eqref{equip3}
under the stronger assumption $(f, g)\in C^\infty_0({\mathbf R}^n)\times 
C^\infty_0({\mathbf R}^n)$
in order to deduce 
\eqref{equip1gen}, \eqref{equip2}
and \eqref{equip3gen} under the weaker assumption 
$(f, g)\in \dot H^1_x\times L^2_x$.

\vspace{0.1cm}

\noindent Since now on we shall assume that
$(f, g)\in C^\infty_0({\mathbf R}^n)\times C^\infty_0({\mathbf R}^n)$.
Following \cite{dg} we introduce the
conformal energy factor
$$[(t^2 + |x|^2) \partial_t u+ 2 n t x_j \partial_j u].$$
Since $\Box u=0$ we get for every $T>0$ the following identity:
$$0=\sum_{j=1}^n 
\int_{0}^T \int [(t^2 + |x|^2) \partial_t u+ 2 t x_j \partial_j u]\Box u 
\hbox{ }dxdt$$
$$=
\int \frac n2 (T^2 + |x|^2) (| \partial_t u(T)|^2+|\nabla_x u(T)|^2)
+2 n T r \partial_{|x|} u(T) \partial_t u(T) \hbox{ } dx 
$$$$- 
\int \frac n2 |x|^2 (|\nabla_x f|^2+|g|^2) \hbox{ } dx
$$$$+ n(n-1) \int_0^T \int t
(|\partial_t u|^2 +|\nabla_x u|^2) \hbox{ } dxdt,$$
where we have used  
the Stokes formula.

\noindent Notice that this identity implies the following inequality:
\begin{equation}\label{stokes3}
\int (T^2 + |x|^2) (| \partial_t u(T)|^2+|\nabla_x u|^2)
+4 T |x| \partial_{|x|} u(T) \partial_t u(T) \hbox{ } dx 
\end{equation}
$$\leq \int |x|^2 (|\nabla_x f|^2+|g|^2) \hbox{ } dx.$$
On the other hand we have the trivial pointwise inequality
$|\partial_{|x|} u|^2\leq |\nabla_x u|^2$ that can be combined with
\eqref{stokes3} in order to give:
\begin{equation}\label{stokes4}
\int (T^2 + |x|^2) (| \partial_t u(T)|^2+|\partial_{|x|} u|^2)
+4 T |x| \partial_{|x|} u(T) \partial_t u(T) \hbox{ } dx 
\end{equation}
$$
\leq \int |x|^2 (|\nabla_x f|^2 + |g|^2) \hbox{ } dx,
$$
and
\begin{equation}\label{stokes5}
\int (T^2 + |x|^2) (| \partial_t u(T)|^2+|\nabla_x u(T)|^2)
-4 T |x| |\nabla_x u(T)||\partial_t u(T)| \hbox{ } dx 
\end{equation}
$$
\leq \int |x|^2 (|\nabla_x f|^2+|g|^2) \hbox{ } dx.
$$
Next recall the basic inequality:
$$(a+b)^2(c+d)^2 + (a-b)^2(c-d)^2
$$
$$\leq 4 (a^2+b^2)(c^2 + d^2) + 16 abcd
\hbox{ } \forall a,b,c,d\in \mathbf R$$
whose proof is completely elementary.
By combining this inequality with
\eqref{stokes4} and \eqref{stokes5} we get
respectively:
$$\int (T+|x|)^2|\partial_t u(T)+ \partial_{|x|}u(T)|^2 + (T-|x|)^2
|\partial_t u(T) -  \partial_{|x|}u(T)|^2 \hbox{ } dx$$$$
\leq 4 \int |x|^2 (|\nabla_x f|^2+ |g|^2) \hbox{ } dx,$$
and
$$\int (T+|x|)^2||\partial_t u(T)| - |\nabla_x u(T)||^2 + (T-|x|)^2
||\partial_t u(T)|+ |\nabla_x u(T)||^2 \hbox{ } dx
$$
$$\leq 4 \int |x|^2 (|\nabla_x f|^2+|g|^2) \hbox{ } dx.$$
This in turn implies:
\begin{equation}\label{equip1fina}\int |\partial_t u(T)+ \partial_{|x|} u(T)|^2 
\hbox{ } dx \leq \frac 4{T^2} \int |x|^2 
(|\nabla_x f|^2+|g|^2) \hbox{ } dx,\end{equation}
\begin{equation}\label{equip2fina}\int ||\partial_t u(T)| 
- |\nabla_x u(T)||^2
\hbox{ } dx 
\leq \frac{4}{T^2} \int |x|^2 (|\nabla_x f|^2+|g|^2) \hbox{ } dx,
\end{equation}
and
\begin{equation}\label{equip3fina}
\int_{2|x|<T} ||\partial_t u(T)| - |\nabla_x u(T)||^2
+ ||\partial_t u(T)| + |\nabla_x u(T)||^2
\hbox{ } dx 
\end{equation}
\begin{equation*}\leq \frac{16}{T^2} \int |x|^2 (|\nabla_x f|^2+|g|^2) \hbox{ } dx.
\end{equation*} 
The proof is complete.

\hfill$\Box$

\end{document}